\documentclass[final]{siamltex}
\usepackage[top=2.95cm,bottom=2.95cm]{geometry}

\usepackage{amsmath} 
\usepackage{amsfonts}
\usepackage{graphicx}
\usepackage{float}
\usepackage{graphicx}
\usepackage{subcaption}
\usepackage{url}
\usepackage{color}
\usepackage{bm}
\usepackage{multirow}

\usepackage{enumitem}

\usepackage{algorithm}
\usepackage{algpseudocode}

\usepackage{tikz}
\usetikzlibrary{trees, cd, babel}

\newcommand{\subscript}[2]{$#1 _ #2$}
\newcommand{\proc}{\text{np}}

\newcommand{\algorithmicbreak}{\textbf{break}}

\DeclareMathOperator*{\argmax}{arg\,max}
\DeclareMathOperator*{\argmin}{arg\,min}

\title{Mixed precision solvers for the all-at-once Runge--Kutta discretization of the heat equation}

\author{Santolo Leveque\footnotemark[1]
        \and Luca Bergamaschi\footnotemark[2]
        \and \'{A}ngeles Mart\'{i}nez\footnotemark[3]
        \and Erin Carson\footnotemark[4]}

\begin{document}

\maketitle
\renewcommand{\thefootnote}{\fnsymbol{footnote}}

\footnotetext[1]{Sokolovsk\'a 49/83, Department of Numerical Mathematics, Faculty of Mathematics and Physics, Charles University, 186 75, Praha 8, Czech Republic (\tt santolo.leveque@matfyz.cuni.cz).}
\footnotetext[2]{Department of Civil Environmental and Architectural Engineering, University of Padova, Via Marzolo 9, 35100, Padova, Italy (\tt luca.bergamaschi@unipd.it).}
\footnotetext[3]{Department of Mathematics, Informatics and Geosciences, University of Trieste, Via Valerio 12/1, 34127, Trieste, Italy (\tt amartinez@units.it).}
\footnotetext[4]{Sokolovsk\'a 49/83, Department of Numerical Mathematics, Faculty of Mathematics and Physics, Charles University, 186 75, Praha 8, Czech Republic {(\tt carson@karlin.mff.cuni.cz)}. \\The first and last authors acknowledge funding from Charles University Research Centre program No. UNCE/24/SCI/005 and from the European Union (ERC, inEXASCALE, 101075632). Views and opinions expressed are those of the authors only and do not necessarily reflect those of the European Union or the European Research Council. Neither the European Union nor the granting authority can be held responsible for them. We acknowledge CINECA for the availability of high-performance computing resources and support under the ISCRA projects IsCe0$\_$P-ROCK and IsCd6$\_$SCALSOL.}

\renewcommand{\thefootnote}{\arabic{footnote}}

\begin{abstract}
We study the effect of mixed precision on the numerical integration of the heat equation discretized with a Runge--Kutta method in time. A full space-time discretization is applied, which results in a very large and sparse linear system to be solved for the numerical approximations and the Runge--Kutta stages of all time steps. The linear system is solved by applying a suitable preconditioned iterative method that can be run in parallel. In order to speed up the solution process, the preconditioner is applied using a mixed precision framework. Sequential results show the robustness of the preconditioner, even when applied in mixed precision. Finally, we present numerical evidence of the improved performance of the mixed precision strategy when applied in a parallel environment, achieving up to a 50\%  reduction in CPU time. 
\end{abstract}

\begin{keywords} 
Time-dependent problems, Parabolic PDEs, Preconditioning, Saddle-point systems, Mixed precision
\end{keywords}

\begin{AMS}
65F08, 65F10, 65N22, 65L06
\end{AMS}

\pagestyle{myheadings}
\thispagestyle{plain}
\markboth{}{MIXED PRECISION FOR RUNGE--KUTTA DISCRETIZATION}

\section{Introduction}\label{sec_1}

Implicit Runge--Kutta methods have been largely studied in the past decades. This class
	of methods is very suited for the long-time integration of very stiff problems
	thanks to their very favourable stability properties. In fact, as opposed to linear
	multistep methods, for which the second Dahlquist barrier states that an A-stable
	linear multistep method cannot have order of convergence greater than two
	\cite[Theorem 6.6]{Lambert}, implicit Runge--Kutta methods of any arbitrary order
	of convergence can be constructed to be A-stable. When solving stiff problems,
	A-stability is an important property a time integration method may possess, as
	this allows a choice of the time step that is completely independent of the problem.
	Another important property that implicit Runge--Kutta methods may have is L-stability, which
	allows the user to employ long time steps without obtaining undesired oscillations
	in the numerical solution. We refer the interested reader to the monographs
	\cite{Butcher, Hairer_Wanner, Lambert} for a thorough exposition on
	Runge--Kutta methods.

The favourable stability properties above make implicit Runge--Kutta methods
	the time integration schemes of choice for stiff problems. For this reason, in recent
	years there has been a renewed interest in employing this class of methods in
	the numerical solution of PDEs. This in turn required researchers to devise fast
	and robust solvers for the large and sparse linear systems arising upon discretization
	with a Runge--Kutta method in time. Simple iterations, such as block Jacobi
	\cite{Mardal_Nilssen_Staff} and block Gauss--Seidel \cite{Staff_Mardal_Nilssen},
	have been employed successfully in the Runge--Kutta integration of the heat equation.
	Strategies based on factorizations of the matrix of the coefficients defining
	the Runge--Kutta method have also been shown to be viable approaches. For example,
	we recall here strategies based on a special LU factorization of the
	coefficient matrix \cite{Rana_Howle_Long_Meek_Milestone}, on the Schur
	decomposition of the coefficient matrix \cite{Southworth_Krzysik_Pazner}, on an SVD
	of the coefficient matrix \cite{Leveque_Bergamaschi_Martinez_Pearson}, and on the
	eigenvalue decomposition of the inverse of the coefficient matrix
	\cite{Southworth_Krzysik_Pazner_Sterck_1}, among many. Further, a
	stage-parallel solver has been devised by employing the diagonalization of the
	lower factor of a special LU factorization \cite{Axelsson_Dravins_Neytcheva};
	numerical results for this strategy in a parallel environment can be found
	in \cite{Munch_Dravins_Kronbichler_Neytcheva}. We also refer to \cite{Durastante_Mazza} for a stage-parallel solver based on the solution of a Sylvester matrix equation with a low-rank correction. Further, we mention monolithic
	multigrid methods as effective and robust strategies for the stage solver; see,
	for example, \cite{AbuLabdeh_MachLachlan_Farrell, Farrell_Kirby_MarchenaMenendez,
	Kirby} and the early works by Vandewalle and co-authors
	\cite{Boonen_VanLent_Vanderwalle, Rosseel_Boonen_Vanderwalle, VanLent_Vanderwalle}.
	Finally, we mention the recent developments of parallel-in-time solvers for the
	all-at-once space-time discretization of the heat equations when employing a
	Runge--Kutta method in time; see \cite{Kressner_Massei_Zhu} and
	\cite{Leveque_Bergamaschi_Martinez_Pearson}. The authors in
	\cite{Kressner_Massei_Zhu} proposed two strategies for the solution of the
	resulting system. First, they observed that the system can be reformulated as
	a Sylvester equation with a right-hand side of low rank. Second, the authors
	observed that the numerical solution can be rewritten
	as a sum of solutions of systems involving $\alpha$-circulant matrices, where
	$\alpha$ is the $j$th root of unity, for some integer $j$. On the other hand,
	the authors in \cite{Leveque_Bergamaschi_Martinez_Pearson} proposed a
	preconditioned iterative method that can be run in parallel thanks to the
	favourable block structure of the preconditioner.

In this work, we consider the numerical integration of the heat equation discretized
	with an implicit Runge--Kutta method in time and finite elements in space. The
	solver is based on the parallel-in-time preconditioning strategy derived
	in \cite{Leveque_Bergamaschi_Martinez_Pearson}, and will be applied in a
	mixed precision framework. The preconditioning strategy devised in
	\cite{Leveque_Bergamaschi_Martinez_Pearson} exploits the block-diagonal structure
	of one of the leading blocks of the discretized system, with the corresponding
	Schur complement resulting in a block-lower triangular matrix; the proposed
	preconditioning strategy is based on an inner stage solver. In what follows, we
	will describe in more detail the preconditioner devised in
	\cite{Leveque_Bergamaschi_Martinez_Pearson}. 

    The emergence of mixed precision computation has become central to modern high-performance computing, enabling substantial gains in performance and energy efficiency while maintaining solution accuracy through carefully designed algorithms and error control strategies; see, e.g., \cite{abdelfattah2021survey, Higham_Mary, kashi2026mixed}. 
    In recent years, researchers have devoted efforts towards devising mixed
	precision Runge--Kutta integrators, see, for example,
	\cite{Balos_Robert_Gardner, Burnett_Gottlieb_Grant, Burnett_Gottlieb_Grant_Heryudono,
	Croci_RosilhodeSouza, Dravins_Koch_Griehl_Kormann, Grant}. Specifically,
	in \cite{Grant} mixed precision is applied to a predictor-corrector Runge--Kutta
	integrator, where the implicit step is solved in lower precision while the explicit
	step is solved in higher precision; numerical evaluation of the stability and of the
	performance of this approach applied to the solution of ODEs has been considered in
	\cite{Burnett_Gottlieb_Grant} and in \cite{Burnett_Gottlieb_Grant_Heryudono},
	respectively, while the work in \cite{Dravins_Koch_Griehl_Kormann} considered performance when the method is applied to PDEs. In \cite{Croci_RosilhodeSouza},
	the authors derive mixed precision Runge--Kutta schemes based on Chebyshev
	polynomials; these methods are explicit integrators, with a stability region that
	grows quadratically with the number of stages. Finally, the authors in
	\cite{Balos_Robert_Gardner} introduced mixed precision in the
	evaluation of exponential integrators. The current work differs from the references
	above, as we are interested in the effect of mixed precision on the linear
	solver. To the best of our knowledge, this is the first study that addresses this
	question.

This work is organized as follows. In Section \ref{sec_2}, we introduce the problem
	and the discretization we employ, then in Section \ref{sec_3} we describe
	the preconditioning strategy adopted in this work. In Section \ref{sec_4} we
	discuss possible mixed precision approaches to be applied within our preconditioning
	strategy. Finally, in Section \ref{sec_5} we report numerical evidence of
	the efficiency of the mixed precision iterative method we propose, drawing
	some conclusions in Section \ref{sec_6}.

\section{The Heat Equation}\label{sec_2}
In what follows, $I_m$ denotes the identity matrix of dimension $m$. Further, we use the following Butcher tableau for defining a Runge--Kutta method:
\begin{displaymath}
\def\arraystretch{1.2}
\begin{array}{c|c}
\mathbf{c}_{\mathrm{RK}} & A_{\mathrm{RK}}\\
\hline
 & \mathbf{b}_{\mathrm{RK}}^\top
\end{array}
\end{displaymath}
where the matrix $A_\mathrm{RK}$ contains the coefficients of the method, and the vectors $\mathbf{b}_{\mathrm{RK}}$ and $\mathbf{c}_{\mathrm{RK}}$ contain, respectively, the weights and the nodes of the integration rule of choice.

Given a spatial domain $\Omega \subset \mathbb{R}^d$, with $d =1, 2, 3$, and a
	final time $T>0$, we consider the numerical integration of the following
	heat equation:
	\begin{equation}\label{heat_equation}
		\left\{
			\begin{array}{rl}
				\vspace{0.25ex}
				\frac{\partial {v}}{\partial t} - \nabla^2 {v} = {f}(\mathbf{x},t)
					& \quad \mathrm{in} \; \Omega \times (0,t_{f}), \\
				\vspace{0.25ex}
				{v}(\mathbf{x},t) = {g}(\mathbf{x},t) & \quad
					\mathrm{on} \; \partial \Omega \times (0,t_{f}),\\
				{v}(\mathbf{x},0) = {v}_0(\mathbf{x}) & \quad
					\mathrm{in} \; \Omega.
			\end{array}
		\right.
	\end{equation}

We employ an $s$-stage Runge--Kutta method in time and then discretize the
	resulting equations with finite elements. Specifically, we discretize the
	time interval $(0,T)$ into $n_t \in \mathbb{N}$ subintervals of length $\tau$.
	Then, with $v_n$ an approximation of $v(\mathbf{x}, t_n)$ with $t_n = n \tau$,
	we evaluate the approximation of the solution at time $t_{n+1}$ as 
	\begin{displaymath}
		v_{n+1} = v_n + \tau \sum_{i=1}^s b_i V_{i,n}.
	\end{displaymath}
	Here, $V_{i,n}$, $i=1, \ldots, s$, are the stages of the method, and we solve
	the system
	\begin{displaymath}
		V_{i,n} - \nabla^2 v_n - \tau \sum_{j=1}^{s} a_{i,j} \nabla^2 V_{j,n}
			= f(\mathbf{x}, t_n + c_i \tau), \quad i=1,\ldots,s.
	\end{displaymath}
	After discretizing with finite elements, we obtain that the numerical solution
	at time $t_n$ is given by
	\begin{equation}\label{update_solution}
		M \mathbf{v}_{n+1} = M \mathbf{v}_n + \tau \sum_{i=1}^s b_i M
			\mathbf{V}_{i,n},
	\end{equation}
	where the discretized stages solve the $s \times s$ block system
	\begin{equation}\label{discretized_system_stages}
		\underbrace{(I_s \otimes M + \tau A_{\mathrm{RK}}\otimes K)}_{\widehat{\Theta}} \mathbf{V}_n = \mathbf{f}_n - K \mathbf{v}_n,
	\end{equation}
	with $\mathbf{V}_n=[\mathbf{V}_{n,1}^\top,\ldots,\mathbf{V}_{n,s}^\top]^\top$,
	and where $\mathbf{f}_n$ contains the discretization of the right-hand side.
	The matrices $M$ and $K$ are the mass and stiffness matrices, respectively.

Starting from a discretization $\mathbf{v}^0$ of the initial condition $v_0$,
	a sequential time-stepping based on a Runge--Kutta scheme solves the system
	in \eqref{discretized_system_stages} and then updates the solution via
	\eqref{update_solution} sequentially for $n=0, \ldots, n_t-1$. By contrast,
	the approach we adopt here is an all-at-once discretization, in which
	\eqref{update_solution} and \eqref{discretized_system_stages} are solved
	simultaneously for all times; that is, we collect \eqref{update_solution}
	and \eqref{discretized_system_stages} for all $n=0,\ldots, n_t-1$, starting
	from the initial condition $\mathbf{v}_0=\mathbf{v}^0$. In matrix form,
	after collecting
	$\mathbf{v} = [\mathbf{v}_{0}^\top,\ldots,\mathbf{v}_{n_t}^\top]^\top$,
	$\mathbf{V} = [\mathbf{V}_{0}^\top,\ldots,\mathbf{V}_{n_t-1}^\top]^\top$, and
	$\mathbf{f} = [\mathbf{f}_{0}^\top,\ldots,\mathbf{f}_{n_t-1}^\top]^\top$, we can write
	\begin{equation}\label{system_all_at_once_RK}
		\underbrace{\left[
			\begin{array}{cc}
				\Phi & \Psi_1\\
				\Psi_2 & \Theta
			\end{array}
		\right]}_{\mathcal{A}}
		\left[
			\begin{array}{c}
				\mathbf{v}\\
				\mathbf{V}
			\end{array}
		\right]=
		\left[
			\begin{array}{c}
				\mathbf{b}\\
				\mathbf{f}
			\end{array}
		\right],
	\end{equation}
	where the vector $\mathbf{b}$ contains information on the initial condition and on
	the boundary conditions. The blocks of the matrix $\mathcal{A}$ are given by
	\begin{equation*}\label{blocks_all_at_once_RK}
		\begin{array}{ll}
			\vspace{1ex}
			\!\!
			\Phi \! = \! \left[
			\!\!
				\begin{array}{cccc}
					M\\
					-M & \ddots\\
					 & \ddots & \ddots\\
					 & & -M & M
				\end{array}
			\!\!
			\right] \!,
			&
			\!\!
			\Psi_1 \! = \! - \! \left[
				\!\!
				\begin{array}{ccc}
					0\\
					\tau \mathbf{b}_{\mathrm{RK}}^\top \otimes M \\
					 & \ddots \\
					 & & \tau \mathbf{b}_{\mathrm{RK}}^\top \otimes M
				\end{array}
				\!\!
			\right] \!,\\
			\!\!
			\Psi_2 \! = \! \left[
				\!\!
				\begin{array}{cccc}
					\mathbf{e} \otimes K\\
					 & \ddots\\
					 & & \mathbf{e} \otimes K & 0
				\end{array}
				\!\!
			\right] \!,
			&
			\!\!
			\Theta = I_{n_t} \otimes \widehat{\Theta} .
		\end{array}
	\end{equation*}
	Here, $\mathbf{e}\in \mathbb{R}^s$ is the column vector of all ones.

\section{Preconditioning Strategy}\label{sec_3}
Given a linear system of the form
	\begin{displaymath}
		\underbrace{\left[
			\begin{array}{cc}
				\Phi & \Psi_1\\
				\Psi_2 & \Theta
			\end{array}
		\right]}_{\mathcal{A}}
		\mathbf{x} = \mathbf{b},
	\end{displaymath}
	with $\Theta$ invertible, we consider as a preconditioner the 
	block-upper triangular matrix
	\begin{equation}\label{preconditioner}
		\mathcal{P}=\left[
			\begin{array}{cc}
				S & \Psi_1\\
				0 & \Theta
			\end{array}
		\right].
	\end{equation}
	Here, $S=\Phi - \Psi_1 \Theta^{-1}\Psi_2$ is the Schur complement.

Preconditioners of the form \eqref{preconditioner} have been studied,
	for instance, in \cite{Ipsen, Murphy:1999:NPI:359189.359190}, and are optimal
	in exact arithmetic. In fact, letting $\sigma(\cdot)$ denote the spectrum of
	a given matrix, under the assumption that $S$ is invertible,
	one can prove that $\sigma(\mathcal{P}^{-1}\mathcal{A})=\left\{ 1 \right\}$
	and the minimal polynomial has degree two. From here, a suitable iterative
	method will converge in at most two iterations.

The preconditioner in \eqref{preconditioner} is not practical, as, excluding
	the effect of inexact arithmetic, the matrix $S$ may be too large to be formed. Further, in some applications the block $\Theta$ may be
	singular, resulting thus in a Schur complement not well defined and a
	preconditioner $\mathcal{P}$ not invertible. For this reason, practical
	preconditioners consist of approximations $\widetilde{\mathcal{P}}$ of
	$\mathcal{P}$, where the diagonal blocks are replaced by suitable matrices.
	Specifically, one considers a block-upper triangular matrix
	\begin{equation}
    \label{AppPrec}
		\widetilde{\mathcal{P}}=\left[
			\begin{array}{cc}
				\widetilde{S} & \Psi_1\\
				0 & \widetilde{\Theta}
			\end{array}
		\right],
	\end{equation}
	where $\widetilde{\Theta}$ is an invertible approximation of
	$\Theta$, and
	$\widetilde{S}\approx \Phi-\Psi_1 \widetilde{\Theta}^{-1}\Psi_2$ is
	an invertible approximation of the Schur complement built with this
	$\widetilde{\Theta}$.

The preconditioner $\widetilde{\mathcal{P}}$ derived
	in \cite{Leveque_Bergamaschi_Martinez_Pearson} for the system
	in \eqref{system_all_at_once_RK} is built upon an approximation
	$\widetilde{\Theta}$ of the inverse of the matrix $\widehat{\Theta}$ defined
	in \eqref{discretized_system_stages}. Specifically, the authors apply a
	fixed number of GMRES \cite{Saad_Schultz} iterations preconditioned with a preconditioner
	$\widetilde{\mathcal{P}}_\mathrm{int}$ based on a SVD of the matrix
	$A_\mathrm{RK}$. Rather than considering a preconditioner based on an SVD of
	$A_\mathrm{RK}$, in our tests we employ a fixed number of GMRES
	iterations preconditioned with the strategy devised in
	\cite{Rana_Howle_Long_Meek_Milestone}. Given this approximation $\widetilde{\Theta}$
	of $\widehat{\Theta}$, the approximation of the
	Schur complement results in the following block-lower triangular matrix:
	\begin{equation}
    \label{AppSchur}
		\widetilde{S}=\left[
			\begin{array}{cccc}
				M\\
				 & \ddots\\
				 &  & M
			\end{array}
		\right]
		\left[
			\begin{array}{cccc}
				I_{n_x}\\
				-I_{n_x}+\widetilde{X} & \ddots\\
				& \ddots & \ddots\\
				& & -I_{n_x}+\widetilde{X} & I_{n_x}
			\end{array}
		\right],
	\end{equation}
	where we set $\widetilde{X}= \tau (\mathbf{b}_{\mathrm{RK}}^\top
	\otimes I_{n_x}) \widetilde{\Theta}^{-1} (\mathbf{e} \otimes K)$.

As a final remark, since the preconditioner $\widetilde{P}$ is non-linear
    (it employs GMRES as inner solver), we have to apply the flexible
    version of GMRES \cite{Saad} as the iterative solver for the system in
    \eqref{system_all_at_once_RK}.

\section{A Mixed Precision Approach}\label{sec_4}

Mixed precision algorithms have attracted significant attention in recent years as an effective response to the difficulty of attaining peak performance on contemporary supercomputers. Building on the classical idea of iterative refinement due to Wilkinson \cite[p. 111]{Wilkinson}, mixed precision arithmetic refers to numerical algorithms that employ different floating-point formats in different stages of the computation to exploit hardware capabilities while retaining high accuracy. In Wilkinson's original scheme, the residual is formed and accumulated in higher precision (for example, double), whereas the correction is obtained by solving a linear system with an LU factorization computed and applied in lower precision (for example, single); this idea has been generalized and analyzed for modern three-precision iterative refinement frameworks by Carson and Higham \cite{carson2018accelerating}. For example, the analysis in \cite{carson2018accelerating} shows that using half precision (fp16) for the expensive LU factorization, double precision for the residual computation, and single precision for all other computations, one can achieve relative accuracy to the level of single precision as long as the infinity-norm condition number of the system $\kappa_\infty(A)<10^8$. Mixed precision implementations of such schemes have demonstrated substantial speedups on modern GPUs and heterogeneous nodes; for dense linear systems on NVIDIA V100 tensor cores, three-precision iterative refinement can deliver up to about $4\times$ faster solves and significantly improved energy efficiency compared with pure double precision methods \cite{haidar2018harnessing, haidar2018design}.

At the system scale, the HPL-MxP benchmark further highlights the growing performance gap between low- and high-precision arithmetic on supercomputers. In the most recent June 2026 TOP500 results, the Aurora supercomputer at Argonne National Laboratory achieved 11.6 Exaflop/s on the mixed precision HPL-MxP benchmark, corresponding to an 11.5$\times$ speedup over the corresponding double precision HPL benchmark \cite{top500-highlights-2026-06}. This widening gap strongly motivates algorithmic designs that carefully combine precisions to improve performance while preserving the desired accuracy. For an overview of recent developments, we refer the interested reader to the surveys \cite{abdelfattah2021survey, Higham_Mary, kashi2026mixed}.

As mentioned above, we consider the integration in a mixed precision framework of
	the heat equation discretized with a Runge--Kutta method in time. We
	are only aware of the works \cite{Balos_Robert_Gardner, Burnett_Gottlieb_Grant,
	Burnett_Gottlieb_Grant_Heryudono, Croci_RosilhodeSouza,
	Dravins_Koch_Griehl_Kormann, Grant} that address Runge--Kutta integration
	in mixed precision arithmetic. As opposed to the references above, in this work
	we investigate the effect of mixed precision in the solution of the system
	in \eqref{system_all_at_once_RK}. To the best of the authors' knowledge, this is
	the first study of mixed precision applied to iterative solvers for the linear
	systems arising from Runge--Kutta discretizations of PDEs. In this regard,
	our framework is related to the studies in \cite{Carson_Dauzickaite, Simoncini_Szyld},
	where the authors analyze the effect of mixed precision on the convergence
	of preconditioned Krylov methods: matrix-vector products and the application of
	the preconditioner have to be applied up to a problem specific precision in order
	to preserve convergence of the method of choice; further, the tolerance on the
	matrix-vector product may be relaxed as the number of iterations proceeds without
	degenerating the attainable accuracy \cite{Simoncini_Szyld}. We leave the study of an iteration-dependent
	matrix-vector tolerance for future work, and rather explore
	the possibility of applying the preconditioner described in Section
	\ref{sec_3} in mixed precision in the tests below. In this regard,
    the analysis on mixed precision flexible GMRES (FGMRES) performed in \cite{Carson_Dauzickaite} gives us some
    guidelines on how to choose the different precisions in the computations. We summarize the main results below.

Suppose we want to solve the linear system $\mathcal{A} \mathbf{x} = \mathbf{b}$ by
    employing preconditioned GMRES or its flexible variant (the latter allowing
    for an iteration-dependent preconditioner) as an iterative solver, with
    $\mathcal{P}_F$ as a preconditioner. In the analysis performed in
    \cite{Carson_Dauzickaite}, the authors have considered the application
    of a split-preconditioned FGMRES method where four precisions are employed
    in the computations. Specifically, if we split the preconditioner
    $\mathcal{P}_F = \mathcal{P}_L \mathcal{P}_R$ into the left preconditioner
    $\mathcal{P}_L$ and the right preconditioner $\mathcal{P}_R$, we can perform computations in four different precisions, that is,
    we have roundoff error $u_{\mathcal{A}}$ for computations with $\mathcal{A}$,
    roundoff error $u_{\mathcal{P}_L}$ for computations with
    $\mathcal{P}_L$, roundoff error $u_{\mathcal{P}_R}$ for computations with
    $\mathcal{P}_R$, and roundoff error $u$ for all other computations.
    Within this framework, a pseudocode for split-preconditioned
    FGMRES is given in Algorithm \ref{4precFGMRES}.

\begin{algorithm}
    \caption{FGMRES with split preconditioning in four precisions}\label{4precFGMRES}
    \begin{algorithmic}
        \State{Given $\mathcal{A} \in \mathbb{R}^{n \times n}$, $\mathbf{b} \in \mathbb{R}^n$, $\mathbf{x}_0 \in \mathbb{R}^n$, left preconditioner $\mathcal{P}_L$, right preconditioner $\mathcal{P}_R$, maximum number of iterations \texttt{maxit}, tolerance \texttt{tol}, precisions $u$, $u_\mathcal{A}$, $u_{\mathcal{P}_L}$, $u_{\mathcal{P}_R}$}
        \State{$\mathbf{t} = \mathcal{A} \mathbf{x}_0$} \Comment{$u_\mathcal{A}$}
        \State{$\mathbf{\tilde{t}} = \mathcal{P}_L^{-1}\mathbf{t}$} \Comment{$u_{\mathcal{P}_L}$}
        \State{$\mathbf{\tilde{b}} = \mathcal{P}_L^{-1}\mathbf{b}$} \Comment{$u_{\mathcal{P}_L}$}
        \State{$\mathbf{r}_0 = \mathbf{\tilde{b}} - \mathbf{\tilde{t}}$} \Comment{$u$}
        \State{$\beta = \| \mathbf{r}_0 \|$, $\mathbf{v}_1 = \mathbf{r}_0 / \beta$} \Comment{$u$}
        \State{$k=0$}
        \While{$k< \mathtt{maxit}$}
            \State{$k = k + 1$}
            \State{$\mathbf{z}_k = \mathcal{P}_R^{-1} \mathbf{v}_k$} \Comment{$u_{\mathcal{P}_R}$}
            \State{$\mathbf{s} = \mathcal{A} \mathbf{z}_k$} \Comment{$u_\mathcal{A}$}
            \State{$\mathbf{w} = \mathcal{P}_L^{-1} \mathbf{s}$} \Comment{$u_{\mathcal{P}_L}$}
            \For{$i=1,\ldots,k$}
                \State{$h_{i,k}= \mathbf{v}_i^\top \mathbf{w}$} \Comment{$u$}
                \State{$\mathbf{w} = \mathbf{w} - h_{i,k} \mathbf{v}_i$} \Comment{$u$}
            \EndFor
            \State{$h_{k+1,k} = \| \mathbf{w} \|$} \Comment{$u$}
            \State{$Z_k = [\mathbf{z}_1, \ldots, \mathbf{z}_k]$, $H_k = \left\{ h_{i,j} \right\}_{1 \leq i \leq j+1,1 \leq j \leq k}$}
            \State{$\mathbf{y}_k = \argmin_\mathbf{y} \| \beta \mathbf{e}_1 - H_k \mathbf{y} \| $} \Comment{$u$}
            \If{$\| \beta \mathbf{e}_1 - H_k \mathbf{y} \| \leq \beta \mathtt{tol}$}
                \State{$\mathbf{x}_k = \mathbf{x}_0 + Z_k \mathbf{y}_k$} \Comment{$u$}
                \State{$t = \mathcal{A} \mathbf{x}_k$} \Comment{$u_\mathcal{A}$}
                \State{$\mathbf{r} = \mathbf{b} - \mathbf{t}$} \Comment{$u$}
                \State{\algorithmicbreak}
            \Else
                \State{$\mathbf{v}_{k+1} = \mathbf{w} / h_{k+1,k}$, $V_{k+1} = [\mathbf{v}_1, \ldots, \mathbf{v}_{k+1} ]$} \Comment{$u$}
            \EndIf
        \EndWhile
    \end{algorithmic}
\end{algorithm}

Bounds on the backward and forward error can be derived, under suitable assumptions, for the four precision split-preconditioned FGMRES method. Specifically, floating point operations involving the inverse of $\mathcal{P}_L$ (resp., of $\mathcal{P}_R$) are assumed to be the sum of the true solution plus a term that can be bounded by the roundoff error $u_{\mathcal{P}_L}$ (resp., $u_{\mathcal{P}_L}$); further, the authors assume that operations involving $\mathcal{P}_L^{-1}\mathcal{A}$ can be written as, excluding higher-order terms, the sum of the true solution plus a vector in norm smaller than a linear combination of the precisions $u_{\mathcal{A}}$ and $u_{\mathcal{P}_L}$. Then, it is possible to prove that the backward error is bounded by a constant (depending only on the dimension of the system and on the current FGMRES iterate) multiplied by a linear combination of the working precision $u$, the precision $u_{\mathcal{A}}$ for the operations with $\mathcal{A}$, and the precision $u_{\mathcal{P}_L}$ for the operations with $\mathcal{P}_L$, with the bound being inversely proportional to a linear combination of the precision $u_{\mathcal{P}_R}$ for the operations with $\mathcal{P}_R$ appearing and the working precision $u$. From this, one can bound the forward error as well, by multiplying the estimate by the condition number of the coefficient matrix $\mathcal{P}_L^{-1}\mathcal{A}$; a bound depending on the split-preconditioned matrix can also be derived, by multiplying the estimate by the condition number of $\mathcal{P}_L^{-1}\mathcal{A}\mathcal{P}_R^{-1}$ times the condition number of $\mathcal{P}_R$, but this results in a weaker bound.

The bounds derived in \cite{Carson_Dauzickaite} can be used not only to decide whether one should apply left-, right-, or split-preconditioning, but they also provide guidance on how to choose the different precisions. For example, if the condition numbers of $\mathcal{A}$ and of $\mathcal{P}_L$ are large, then one may require $u_{\mathcal{A}} \ll u$; further, for ill-conditioned $\mathcal{A}$ and $\mathcal{P}_L$, in order to have a small backward error one may require $u_{\mathcal{A}} \ll u_{\mathcal{P}_L}$; notably, the precision $u_{\mathcal{P}_R}$ for the operations with $\mathcal{P}_R$ does not influence the backward error, as long as it is bounded by the norm of the inverse of $\mathcal{P}_R$ divided by the norm of the error matrix resulting from matrix-vector products with $\mathcal{P}_R$ (for ill-conditioned $\mathcal{P}_R$ a small value of $u_{\mathcal{P}_R}$ may be required for this condition to be satisfied).

Regarding the preconditioning approach described in Section \ref{sec_3}, when applied in exact arithmetic, the preconditioned matrix $\mathcal{P}^{-1}\mathcal{A}$ is block-Jordan. An analysis of the condition number of a matrix of this form is beyond the scope of this work. The analysis becomes even more involved if the preconditioner is applied inexactly, as we do here. Rather, we are interested in the effect of different precisions within the application of the preconditioned iterative method. Specifically, in our experiments we study the behavior of the preconditioner $\widetilde{\mathcal{P}}$ defined in \eqref{AppPrec} when applied in a different precision than the working precision, either within right-preconditioned FGMRES or within left-preconditioned FGMRES.

As a part of the application of the preconditioner, one
	has to apply an inexact inverse of matrices of the form $M+c K$, with $c$ a given
	constant; this will be done by means of a suitable multigrid routine.
	Recent studies on the mixed precision multigrid
	solvers can be found in \cite{McCormick_Benzaken_Tamstorf,
	Tamstorf_Benzaken_McCormick, Vacek_Anzt_Carson}. In \cite{McCormick_Benzaken_Tamstorf,
	Tamstorf_Benzaken_McCormick}, the authors analyze the effect of the use of
	different precision at each grid of the V-cycle, decreasing the precision for
	coarser levels, presenting error bounds for the framework considered.
	In \cite{Vacek_Anzt_Carson}, an analysis of the application in mixed precision
	of the smoother and of the coarse solver of a V-cycle is performed; theoretical
	bounds on the error show the effect of the precision of the smoother and of
	the coarse solver in the convergence of the method. Overall, by employing suitable
	precisions for each level of the grid and for the smoother and the coarse solve
	one can still obtain a convergent V-cycle. Further, numerical results show the
	speed-up achieved by employing a mixed precision multigrid framework, see \cite{Vacek_Anzt_Carson, Vacek_Carson_Soodhalter}. Although
	these are viable approaches, in the numerical tests below we only consider
	a multigrid applied in a given precision.

\section{Numerical Results}\label{sec_5}
We now present numerical results of our approach when solving the heat
	equation discretized with Radau IIA methods in time. We compare the
	preconditioned iterative method described in Section \ref{sec_3},
	run in double precision, with
	the same preconditioner applied in a mixed precision framework.
	Specifically, one can consider the following cases:
	\begin{enumerate}[label=(\subscript{A}{{\arabic*}})]
		\item the system in \eqref{system_all_at_once_RK} is solved in double
			precision, with the preconditioners $\widetilde{\mathcal{P}}$ and
			$\widetilde{\mathcal{P}}_\mathrm{int}$ applied in double
			precision;\label{approach_A1}

		\item the system in \eqref{system_all_at_once_RK} is solved in
			double precision, with the preconditioners $\widetilde{\mathcal{P}}$
			and $\widetilde{\mathcal{P}}_\mathrm{int}$ applied in single
			precision;\label{approach_A2}

		\item the system in \eqref{system_all_at_once_RK} is solved in
			single precision, with the preconditioners $\widetilde{\mathcal{P}}$
			and $\widetilde{\mathcal{P}}_\mathrm{int}$ applied in single
			precision;\label{approach_A3}

		\item the system in \eqref{system_all_at_once_RK} is solved in
			double precision, with the preconditioners $\widetilde{\mathcal{P}}$
			and $\widetilde{\mathcal{P}}_\mathrm{int}$ applied in half
			precision;\label{approach_A4}

		\item the system in \eqref{system_all_at_once_RK} is solved in
			double precision, with the preconditioner $\widetilde{\mathcal{P}}$
			applied in single precision, while the preconditioner
			$\widetilde{\mathcal{P}}_\mathrm{int}$ is applied in half
			precision.\label{approach_A5}
	\end{enumerate}

From numerical evidence, we observed that approach \ref{approach_A3} results
	in a loss of order of convergence of the method, when employing
	higher-order spatial discretizations; this is not surprising, as the use
	of single precision for storing the linear system results in a loss
	of information. Further, within the approach \ref{approach_A4} the
	preconditioner loses its robustness, resulting in a number of iterations
	dependent on the discretization and on the number of stages of the
	Runge--Kutta method; this is in accordance with the findings in
	\cite{Carson_Dauzickaite}, as the precision for computations with
	the preconditioner has to be chosen carefully in order to preserve the
	quality of the preconditioner. For these reasons, we will not consider
	these two approaches in what follows.

Regarding approach \ref{approach_A5}, we would like to mention that this
	requires the solver to be run on modern architectures that provide a half
	precision framework, for example, on GPUs. Unfortunately, the multigrid
	employed in our solver is not implemented on GPUs. For this reason, we
	present numerical results for \ref{approach_A5} only in a sequential
	framework, just as a proof of concept of this viable approach.

Given the above, in what follows we consider approach \ref{approach_A1},
	\ref{approach_A2}, and \ref{approach_A5} in a sequential framework. Then,
	we compare approach \ref{approach_A1} with approach \ref{approach_A2} within
	a parallel environment.

We use $\mathbf{Q}_1$ and $\mathbf{Q}_2$ elements as spatial discretization.
	For each level of refinement $l$, we set as mesh-size $h = 2^{1-l}$ for $Q_1$
	elements and $h = 2^{-l}$ for $Q_2$ elements. Regarding the time grid, letting
	$q_{\mathrm{FE}}$ be the order of the finite element approximations and
	$q_{\mathrm{RK}}$ be the order of the Runge--Kutta method we employ, we set
	$n_t$ to be the closest integer such that
	$\tau \leq h^{q_{\mathrm{FE}}/q_{\mathrm{RK}}}$, where $\tau=\frac{t_f}{n_t}$ is
	the time-step.

\begin{figure}[tbhp]
	\footnotesize
	\centering
	\begin{tikzpicture}[%
		every node/.style={draw=black, thick, anchor=west},
		grow via three points={one child at (-0.0,-0.7) and two children at (-0.0,-0.7) and (-0.0,-1.4)},
		edge from parent path={(\tikzparentnode.210) |- (\tikzchildnode.west)}]
		\node {Krylov solver: FGMRES}
		child {node at (-0.2,-0.15) [align=left] {Solve for
			$\mathcal{A}$ with block triangular\\
			preconditioner $\widetilde{\mathcal{P}}$ defined in Section \ref{sec_3}}
			child {node at (-0.2,-0.15) [align=left] {Solve for $\widetilde{\Theta}$: block-diagonal solve, approximating inverse of $\widehat{\Theta}$}
				child {node at (-0.2,-0.1) [align=left] {Krylov solver: GMRES}
				child {node at (-0.2,-0.1) [align=left]{Apply Rana's preconditioner}}
                child {node at (-0.2,-0.2) [align=left]{Apply V-cycle for solving system\\
                of the form $M+a K$}}
			}}
			child {node at (-0.2,-2.8) {Solve for Schur complement approximation $\widetilde{S}$}
				child { node at (-0.2,-0.1) [align=left] {Block-diagonal solve,\\
                approximating inverse of $M$}
				child {node at (-0.2,-0.2) [align=left] {Chebyshev semi-iteration for\\
                mass matrix}}}
				child {node at (-0.2,-1.2) {Block-forward substitution, approximating inverse of $\widehat{\Theta}$}
					child {node at (-0.2,-0.1) [align=left] {Krylov solver: GMRES}
				child {node at (-0.2,-0.1) [align=left]{Apply Rana's\\ preconditioner}}
                child {node at (-0.2,-0.4) [align=left]{Apply V-cycle 
                for\\ solving
                system of\\ the form $M+a K$}}
			}
					child[missing]{}
					child[missing]{}
				child[missing]{}
				}
			}
			child[missing]{}
			child[missing]{}
			child[missing]{}
		};
	\end{tikzpicture}
	\caption{Solver diagram for the preconditioner $\widetilde{\mathcal{P}}$.}
	\label{fig:prec}
\end{figure}
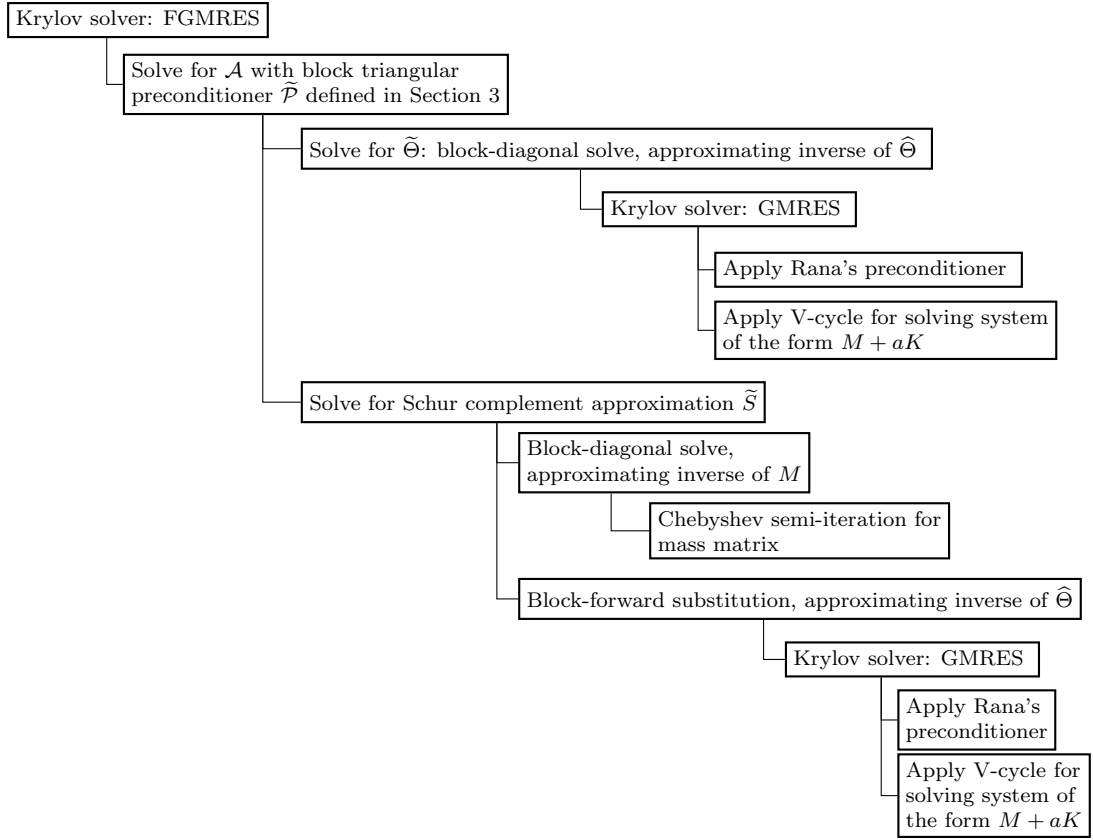

We employ FGMRES restarted every 20 iterations as an outer solver, seeking
	a reduction of $10^{-8}$ on the relative residual and allowing for a maximum
	of 100 iterations. As we mentioned above, the system in
	\eqref{discretized_system_stages} is solved by employing a fixed number
	of GMRES iterations (3 in our tests) preconditioned by Rana's preconditioner \cite{Rana_Howle_Long_Meek_Milestone}.
	Each application of the outer preconditioner requires the approximate inversion of
	a mass matrix, which is done by applying 20 steps of Chebyshev
	semi-iteration \cite{GolubVargaI,GolubVargaII,Wathen_Rees} with
	Jacobi splitting, and the approximate inversion of the diagonal blocks
	within Rana's preconditioner, which is done by applying 3 cycles of
	the \texttt{HSL\_MI20} solver \cite{HSL_MI20}. We summarize the application of the preconditioner in a sequential framework in Figure \ref{fig:prec}. We mention that in the parallel implementation of the preconditioner, the block-forward substitution within the Schur complement approximation $\widetilde{S}$ is replaced by a suitable MGRIT routine.

\subsection{Sequential Results}
We begin by reporting results when running the solver sequentially on a laptop.
	All tests are run on MATLAB R2018b, using a 1.70GHz Intel quad-core i5
	processor and 8 GB RAM on an Ubuntu 18.04.1 LTS operating system.
    We recall that in this section we compare approaches
    \ref{approach_A1}, \ref{approach_A2}, and \ref{approach_A5}. Here, we apply right-preconditioned FGMRES as the outer solver.

We solve problem \eqref{heat_equation} with exact solution given by
    \begin{displaymath}
        v(\boldsymbol x, t) = \exp\left(t_f -t\right) \cos\left(\frac{\pi x_1}{2}
            \right)\left(\frac{\pi x_2}{2} \right) + 1;
    \end{displaymath}
    initial condition $v_0$, boundary conditions $g$, and force function
    $f$ are derived by this choice of $v$. Here, we integrate the problem up to
    $t_f=2$. Below, the discretization error
    $v_{\text{error}}$ is defined as
	\begin{displaymath}
		v_{\text{error}} = \max_{n} \left\{ \dfrac{| v_{\mathtt{j},n} -
			v_{\mathtt{j},n}^{\mathrm{\: sol}}
			|}{|v_{\mathtt{j},n}^{\mathrm{\: sol}}|}, \; \mathrm{with} \;
			\mathtt{j} = \argmax_j | v_{j,n} - v_{j,n}^{\mathrm{\: sol}} | \right\},
	\end{displaymath}
	with $v_{j,n}$ and $v_{j,n}^{\mathrm{\: sol}}$ being the entries of the
	numerical solution $\mathbf{v}$ and the (discretized) exact solution for $v$
	at time $t=t_n$.

We report the degrees of freedom (DoF), the number of FGMRES iterations,
	and the relative errors of the discretization in Table
	\ref{table_sequential_double}, when solving
	the linear system in double precision and applying the preconditioner in
	double precision (approach \ref{approach_A1}), while we report the number
	of FGMRES iterations and the relative errors of the
	discretization in Table
	\ref{table_sequential_single}, when solving
	the linear system in double precision and applying the preconditioner in
	single precision (approach \ref{approach_A2}). Further, we report the number
	of FGMRES iterations and the relative errors of
	the discretization in Table
	\ref{table_sequential_mixed}, when applying
	the preconditioner on 3 levels of precision (approach \ref{approach_A5}).
	Single and half precision are simulated by employing the chop algorithm
	\cite{Higham_Pranesh}.

\begin{table}[!ht]
\caption{All-at-once solve of the heat equation: degrees of freedom (DoF),
	FGMRES iterations, and resulting relative errors in $v$ for Radau
	IIA methods, employing $\mathbf{Q}_1$ and $\mathbf{Q}_2$ finite elements,
	when applying the preconditioner in double
	precision.}\label{table_sequential_double}
	\begin{center}
		\begin{footnotesize}
			{\begin{tabular}{c|c|ccc|ccc}
				\multicolumn{1}{c}{} & \multicolumn{1}{c|}{} &
					\multicolumn{3}{c|}{$\mathbf{Q}_1$} &
					\multicolumn{3}{c}{$\mathbf{Q}_2$} \\
				\cline{3-8}
				$s$ & $l$ & DoF & $\texttt{it}$ & $v_{\text{error}}$ & DoF &
					$\texttt{it}$ & $v_{\text{error}}$ \\
				\hline
				\multirow{4}{*}{2} & $3$ & 931        & 4    & 5.09e-03
                				         & 5625       & 5     & 2.02e-04 \\
				                   & $4$ & 5625       & 5     & 1.35e-03
				                         & 47,089     & 5     & 2.51e-05 \\
				                   & $5$ & 38,440     & 5     & 3.46e-04
				                         & 384,993    & 5     & 3.28e-06 \\
				                   & $6$ & 254,016    & 5     & 8.35e-05
				                         & 3,112,897  & 5      & 4.03e-07 \\
				                   & $7$ & 1,564,513  & 5     & 2.05e-05
				                         & 25,034,625 & 5     & 6.00e-08 \\
				\hline
				\multirow{4}{*}{3} & $3$ & 833       & 6     & 5.71e-03
				                         & 4725      & 7     & 2.57e-05 \\
				                   & $4$ & 4725      & 7     & 1.46e-03
				                         & 27,869    & 7     & 2.27e-06 \\
				                   & $5$ & 27,869    & 7     & 3.30e-04
				                         & 178,605   & 7     & 2.07e-07 \\
				                   & $6$ & 130,977   & 7     & 1.03e-04
				                         & 1,048,385 & 7      & 3.05e-08 \\
				                   & $7$ & 725,805   & 6     & 2.36e-05
				                         & 6,567,525 & 7      & 3.02e-08 \\
				\hline
				\multirow{4}{*}{4} & $3$ & 784       & 8     & 5.59e-03
				                         & 4725      & 9     & 1.94e-05 \\
				                   & $4$ & 4725      & 8     & 1.54e-03
				                         & 24,986    & 10     & 1.17e-06 \\
				                   & $5$ & 24,986    & 9     & 3.77e-04
				                         & 142,884   & 10     & 7.19e-08 \\
				                   & $6$ & 123,039   & 9     & 9.00e-05
				                         & 741,934   & 10      & 7.77e-09 \\
				                   & $7$ & 580,644   & 9     & 2.12e-05
				                         & 3,966,525 & 10      & 6.47e-09 \\
				\hline
				\multirow{4}{*}{5} & $3$ & 931       & 7     & 5.91e-03
				                         & 5625      & 7     & 1.92e-05 \\
				                   & $4$ & 5625      & 7     & 1.54e-03
				                         & 24,025    & 8     & 1.20e-06 \\
				                   & $5$ & 24,025    & 7     & 4.01e-04
				                         & 146,853   & 8     & 7.16e-08 \\
				                   & $6$ & 123,039   & 7     & 9.65e-05
				                         & 693,547   & 8      & 2.50e-09 \\
				                   & $7$ & 596,773   & 8     & 2.27e-05
				                         & 3,186,225 & 8      & 1.04e-08 
			\end{tabular}}
		\end{footnotesize}
	\end{center}
\end{table}

\begin{table}[!ht]
\caption{All-at-once solve of the heat equation:
	FGMRES iterations and resulting relative errors in $v$ for Radau
	IIA methods, employing $\mathbf{Q}_1$ and $\mathbf{Q}_2$ finite elements,
	when applying the preconditioner in single
	precision.}\label{table_sequential_single}
	\begin{center}
		\begin{footnotesize}
			{\begin{tabular}{c|c|cc|cc}
				\multicolumn{1}{c}{} & \multicolumn{1}{c|}{} &
					\multicolumn{2}{c|}{$\mathbf{Q}_1$} &
					\multicolumn{2}{c}{$\mathbf{Q}_2$} \\
				\cline{3-6}
				$s$ & $l$ & $\texttt{it}$ & $v_{\text{error}}$ & $\texttt{it}$ &
					$v_{\text{error}}$ \\
				\hline
				\multirow{4}{*}{2} & $3$ & 4      & 5.09e-03
				                         & 5      & 2.02e-04 \\
				                   & $4$ & 5      & 1.35e-03
				                         & 5      & 2.51e-05 \\
				                   & $5$ & 5      & 3.46e-04
				                         & 5      & 3.29e-06 \\
				                   & $6$ & 5      & 8.35e-05
				                         & 5       & 4.45e-07 \\
				                   & $7$ & 5      & 2.05e-05
				                         & 5      & 4.43e-08 \\
				\hline
				\multirow{4}{*}{3} & $3$ & 6      & 5.71e-03
				                         & 7      & 2.57e-05 \\
				                   & $4$ & 7      & 1.46e-03
				                         & 7      & 2.27e-06 \\
				                   & $5$ & 7      & 3.30e-04
				                         & 7      & 2.07e-07 \\
				                   & $6$ & 7      & 1.03e-04
				                         & 7       & 3.03e-08 \\
				                   & $7$ & 6      & 2.36e-05
				                         & 7       & 2.34e-08 \\
				\hline
				\multirow{4}{*}{4} & $3$ & 8      & 5.59e-03
				                         & 9      & 1.94e-05 \\
				                   & $4$ & 8      & 1.54e-03
				                         & 10      & 1.17e-06 \\
				                   & $5$ & 9      & 3.77e-04
				                         & 10      & 7.19e-08 \\
				                   & $6$ & 9      & 9.00e-05
				                         & 10       & 7.78e-09 \\
				                   & $7$ & 9      & 2.12e-05
				                         & 10       & 6.31e-09 \\
				\hline
				\multirow{4}{*}{5} & $3$ & 7      & 5.91e-03
				                         & 7      & 1.92e-05 \\
				                   & $4$ & 7      & 1.54e-03
				                         & 8      & 1.20e-06 \\
				                   & $5$ & 7      & 4.01e-04
				                         & 8      & 7.17e-08 \\
				                   & $6$ & 7      & 9.65e-05
				                         & 8       & 2.54e-09 \\
				                   & $7$ & 8      & 2.27e-05
				                         &        & 1.14e-08 
			\end{tabular}}
		\end{footnotesize}
	\end{center}
\end{table}

\begin{table}[!ht]
\caption{All-at-once solve of the heat equation:
	FGMRES iterations and resulting relative errors in $v$ for Radau
	IIA methods, employing $\mathbf{Q}_1$ and $\mathbf{Q}_2$ finite elements,
	when applying the preconditioner in mixed
	precision.}\label{table_sequential_mixed}
	\begin{center}
		\begin{footnotesize}
			{\begin{tabular}{c|c|cc|cc}
				\multicolumn{1}{c}{} & \multicolumn{1}{c|}{} &
					\multicolumn{2}{c|}{$\mathbf{Q}_1$} &
					\multicolumn{2}{c}{$\mathbf{Q}_2$} \\
				\cline{3-6}
				$s$ & $l$ & $\texttt{it}$ & $v_{\text{error}}$ & $\texttt{it}$
					& $v_{\text{error}}$ \\
				\hline
				\multirow{4}{*}{2} & $3$ & 6      & 5.09e-03
				                         & 6      & 2.02e-04 \\
				                   & $4$ & 6      & 1.35e-03
				                         & 6      & 2.51e-05 \\
				                   & $5$ & 6      & 3.46e-04
				                         & 7      & 3.29e-06 \\
				                   & $6$ & 6      & 8.35e-05
				                         & 9       & 4.34e-07 \\
				                   & $7$ & 7       & 2.06e-05
				                         & 17      & 8.36e-08 \\
				\hline
				\multirow{4}{*}{3} & $3$ & 6      & 5.71e-03
				                         & 7      & 2.57e-05 \\
				                   & $4$ & 7      & 1.46e-03
				                         & 7      & 2.28e-06 \\
				                   & $5$ & 7      & 3.30e-04
				                         & 7      & 2.08e-07 \\
				                   & $6$ & 7      & 1.03e-04
				                         & 8       & 1.84e-08 \\
				                   & $7$ & 8      & 2.37e-05
				                         & 13      & 4.50e-08 \\
				\hline
				\multirow{4}{*}{4} & $3$ & 8      & 5.59e-03
				                         & 9      & 1.94e-05 \\
				                   & $4$ & 9      & 1.54e-03
				                         & 10      & 1.17e-06 \\
				                   & $5$ & 9      & 3.77e-04
				                         & 10      & 6.97e-08 \\
				                   & $6$ & 10      & 9.00e-05
				                         & 10       & 4.08e-09 \\
				                   & $7$ & 10      & 2.12e-05
				                         & 13      & 1.02e-08 \\
				\hline
				\multirow{4}{*}{5} & $3$ & 7      & 5.91e-03
				                         & 8      & 1.92e-05 \\
				                   & $4$ & 8      & 1.54e-03
				                         & 8      & 1.20e-06 \\
				                   & $5$ & 8      & 4.01e-04
				                         & 8      & 7.24e-08 \\
				                   & $6$ & 8      & 9.65e-05
				                         & 9       & 3.60e-09 \\
				                   & $7$ & 9      & 2.27e-05
				                         & 12       & 7.22e-08 
			\end{tabular}}
		\end{footnotesize}
	\end{center}
\end{table}

From Tables \ref{table_sequential_double}--\ref{table_sequential_mixed}, we can observe
	the robustness of the proposed strategy. The double precision solver can
	achieve convergence in at most 10 iterations for the problem considered here,
	independently of the discretization parameters and of the number of stages of
	the Runge--Kutta method. Interestingly, the convergence of the solver is not
	affected when applying the preconditioner in single precision, with the
	discretization error being polluted only mildly for $\mathbf{Q}_2$ elements;
    this may be because the chop algorithm is handling numbers as double and
    not as single. As we will see in the following section, an increase in the
    number of outer FGMRES iterations may be expected when applying the
    preconditioner in single precision. Regarding approach \ref{approach_A5}, the solver shows some dependence on the level of
    refinement adopted, in particular when applying higher order finite elements,
    see the right column in Table
	\ref{table_sequential_single}. Further, the discretization error is more polluted
	in this case by the mixed precision preconditioning approach, but it remains
	comparable to the one obtained by the fully double precision approach.
	Regarding the discretization
	error, we observe the predicted order of convergence of the method, until it is
	slightly polluted by the linear tolerance of the solver for the finest grid of
	$\mathbf{Q}_2$ elements, see, for example, levels $l=6$ and $l=7$ when employing
	5 stages in Table \ref{table_sequential_mixed}. 

Given the promising sequential results shown here, in the following
    section we will study the speed-up obtained on parallel
    architectures when considering a mixed precision approach compared to
    the full double precision solver.

\subsection{Parallel Results}
We now consider a parallel implementation of the proposed mixed precision preconditioner. In what follows, we apply left-preconditioned FGMRES as the outer solver.

The selected test case regards the heat equation, following
    \cite[Section 4.1.1]{Leveque_Bergamaschi_Martinez_Pearson},
    in which we slightly changed the exact solution to handle large final times, as follows:
    \[ v(\boldsymbol x, t) = \exp\left(\frac{t_f -t}{10}\right) \cos\left(
        \frac{\pi x_1}{2} \right)\left(\frac{\pi x_2}{2} \right) + 1.\]
    The relevant parameters are summarized in Table \ref{tab1}.
    \begin{table}[h!]
        \caption{Parameters of the simulation}\label{tab1}
        \begin{center}
            \begin{tabular}{l|l|l}
            	\multicolumn{3}{c}{ outer FGMRES  \texttt{tol} = $10^{-8}$ \ \texttt{restart} $= 15$ }\\
                \hline
            	& Test case \#1 & Test case \#2 \\
            	\hline
            	$\Delta t$  & $8.236 \times 10^{-2}$ &	$5.437 \times 10^{-2}$ \\
            	$h$ &   $1.563 \times 10^{-2}$ & 	 $7.813 \times 10^{-3}$ \\
            	Tfin &  21.084 & 27.838 \\
            	n. of stages $s$ &            3 &  3 \\
                RK method &            Radau IIA &  Radau IIA \\
                RK order &            5 & 5 \\
                FE space & $Q_2$  (order 3 accuracy) & $Q_2$   \\
            	$n_x$ &        65025 &  261121 \\
            	$n_t$ &          256 &  512 \\
        	    All-at-once DOF &     66\,650\,624 & 535\,036\,928 
            \end{tabular}
    	\end{center}
    \end{table}

We run the MPI-Fortran90 all-at-once solver on the Leonardo Booster module
    (located at Cineca, Bologna, Italy), which features 3,456 nodes, each with one
    32-core Intel Xeon Platinum 8358 CPU (Ice Lake), 512 GB RAM, and 4x NVIDIA Ampere A100
    64GB HBM2e GPUs; the solver is equipped with the preconditioner for the stages as in
    \cite{Rana_Howle_Long_Meek_Milestone},  using either double or single precision
    \emph{for the preconditioner application only}, i.e. the solution of a linear
    system having $\widetilde P =
    \begin{bmatrix} \widetilde S & \Psi_1 \\ 0 & \widetilde \Theta \end{bmatrix}$,
    see \eqref{AppPrec}, as the coefficient matrix, by using $\proc = 1, 8, 64$ processors.

Notice that the outer (F)GMRES solver requires the solution of several linear systems with
    the matrix $\widetilde \Theta$, which in turn is obtained by an inner GMRES solver. We
    had to tune both the maximum number of iterations and the tolerance for this inner solver,
    depending on the precision. In particular, we set:
    \begin{itemize}
    	\item (Single precision) \texttt{maxit} $ = 3$, \texttt{tol} = $10^{-3}$.
    	\item (Double precision) \texttt{maxit} $ = 5$, \texttt{tol} = $10^{-6}$.
    \end{itemize}
    In Tables \ref{tab:case1} and \ref{tab:case2}, we report the inner/outer number
    of iterations and the CPU times for the most time-consuming operations.  We can notice
    the outstanding parallelization degree in the solution of the $(2,2)$ block
    $\widetilde \Theta$, whereas the solution of the Schur complement matrix conveys
    the solution of the bidiagonal block system in \eqref{AppSchur}, yet
    efficiently parallelized by the Xbraid solver \cite{XBraid}.  We also report the
    parallel speedup computed as 
    \[ S_{\proc}  = \frac{T_{\proc}}{T_1},\]
    being $T_{\proc}$ the total elapsed CPU time on $\proc$ processors.

\begin{figure}[h!]
	\vspace{-2.4cm}
    \begin{minipage}{7.5cm}
    	\includegraphics[width=7.5cm]{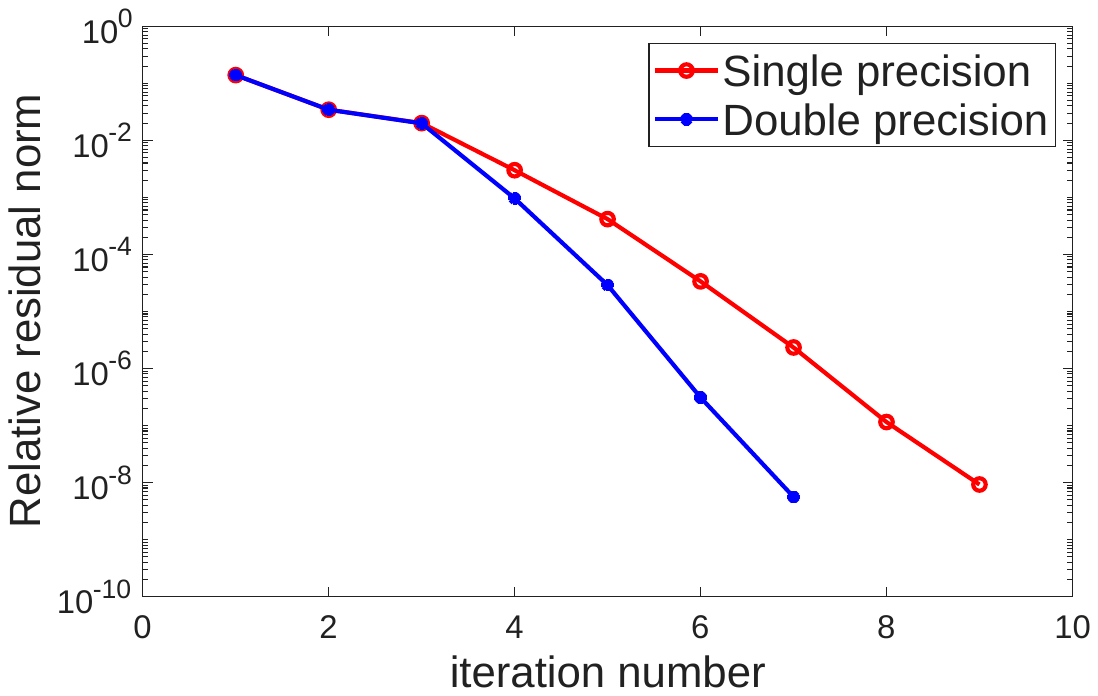}
    \end{minipage}
	\hspace{-4mm}
    \begin{minipage}{7.5cm}
    	\includegraphics[width=7.5cm]{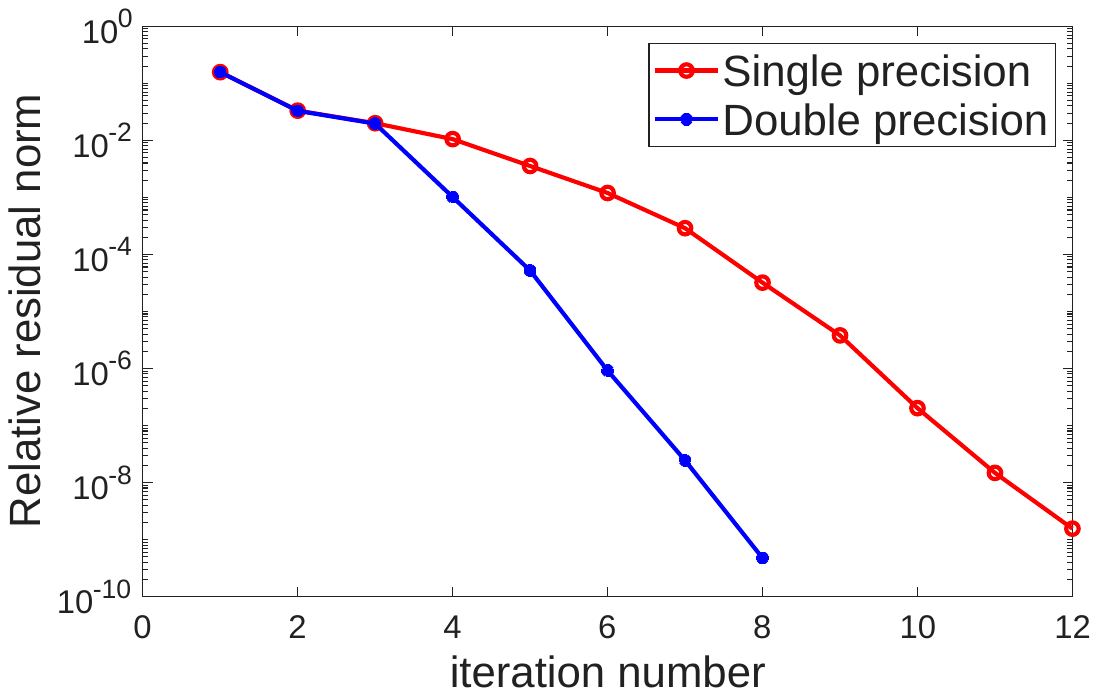}
    \end{minipage}
	\vspace{-2.8cm}
	\caption{Convergence profiles. On the left, test case \#1, on the right test case \#2.}
	\label{convprof}
\end{figure}

The different accuracy in the preconditioner (depending on the different tolerances for
    the inner GMRES solver applied to matrix $\widetilde \Theta$) moderately impacts the
    number of outer iterations, which increase by 30\% -- 50\%, depending
    on the test case, with slightly different convergence profiles as shown in
    Figure \ref{convprof}. However, the CPU times are always significantly
    lower using the single-precision preconditioner, irrespective of the number of
    processors employed, see also Table \ref{Single_vs_Double}.

\begin{table}[h!]
	\caption{Test case \#1. Timings and speedups of the parallel simulations.}
	\label{tab:case1}
    \begin{center}
    	\begin{tabular}{l|r|rrr|r|rr}
		      &	  & \multicolumn{3}{c}{CPU}  &$S_{\proc}$ &\multicolumn{1}{c}{its} \\
    		&	\proc & $\widetilde \Theta $ & $\widetilde S$ & $T_{\proc}$ &  & (inner)\\
	       	\hline
    		&	1 &1085.0 & 1349.6 & 2469.4 &-- &9 (24) \\
    		single		&	8 & 184.9 & 273.6  &  464.9 &5.3& 9 (24) \\
    		&64 & 24.5 & 107.4  &  135.0  & 18.3&9 (24)\\
            \hline
    		&1& 1446.0 & 1883.4& 3316.0 & -- &7  (23) \\
    		double          & 8  &225.7&335.8& 585.3 & 5.7 & 7  (23)\\
    		&64&  30.2  & 119.1& 149.9  &22.1 &7 (23)
        \end{tabular}
    \end{center}
\end{table}

\begin{table}[h!]
	\caption{Test case \#2. Timings and speedups of the parallel simulations.}
	\label{tab:case2}
    \begin{center}
    	\begin{tabular}{l|r|rrr|r|rr}
    		&	  & \multicolumn{3}{c}{CPU}  &$S_{\proc}$ &\multicolumn{1}{c}{its} \\
    		&	\proc & $\widetilde \Theta $ & $\widetilde S$ & $T_{\proc}$ &  & (inner)\\
    		\hline
    		&	1 & 18423.3& 17418.1 &  36272.4&-- &11 (33)\\
    		single		&	8 &2840.6  & 3765.2  & 6683.2 & 5.4&12 (34) \\
    		&64 & 387.7  & 1247.9  &  1677.8  & 21.6&12 (34)\\
            \hline
    		&1& 34258.0 & 43328.8  &  77926.6& & 8 (30) \\
    		double         & 8 &5167.7& 6535.5 &11745.0 &6.6 & 8 (30)\\
    		&64&  724.2  & 1861.2 & 2592.2 & 30.1& 8 (30)
        \end{tabular}
    \end{center}
\end{table}

We would like to mention that the speed-up obtained in single precision is not influenced by the number of the maximum number of iterations nor by the tolerance adopted for the inner solver. To show this, we run test case \#2, with the preconditioner applied in double precision, employing the same parameters used in the single precision runs (that is, \texttt{maxit}$ = 3$ and \texttt{tol} = $10^{-3}$). The results, with 8 processors, are shown in the first and second rows of Table \ref{tab:np8}, whereas in the second row of Table \ref{tab:np8} we also report the previous results with $\proc = 8$ processors, for the sake of comparison. We note that the parallel speedups are slightly larger using double instead of single precision. This is because floating-point operations in double precision are roughly twice as expensive as floating-point operations in single, but the reduction in communication time is not as significant.

\begin{table}[h!]
    \caption{Effects of lower tolerance for the inner solver
        and reduced precision with \proc = 8.}
    \label{tab:np8}
	\begin{center}
        \begin{tabular}{l|r|rrr|r|rr}
            &	  & \multicolumn{3}{c}{CPU}  &\multicolumn{1}{c}{its} & & \\
            &	\proc & $\widetilde \Theta $ & $\widetilde S$ & $T_{\proc}$ & (inner)& \texttt{maxit} & \texttt{tol} \\
            single	&	8 &2840.6& 3765.2 & 6683.2 &12 (34)& $3$ & $10^{-3}$  \\
            double  &       8 &3922.3& 5277.6 & 9257.0 &11 (27)& $3$ & $10^{-3}$  \\
            double  &       8 &5167.7& 6535.5 &11745.0 & 8 (30)& $5$ & $10^{-6}$ 
        \end{tabular}
	\end{center}
\end{table}

For the last numerical test (denoted with test case \#3) we employed a 5-stage Radau IIA method in combination with
    Q1 Finite Elements in space with a very fine spatial discretization. The parameters are
    now: $\Delta t = 0.2912$, Tfin = $149.1, h = 3.91 \times 10^{-3}, n_x = 261121$
    and $n_t = 512$. Overall, the number of DOF is  802\,424\,832. The GMRES parameters are
    the same as in the previous examples,  with the only exception of the inner
    GMRES \texttt{maxit} which has been set to 4, to follow the (yet slight) dependence on
    the number of stages of the inner preconditioner we employed, following
    \cite{Rana_Howle_Long_Meek_Milestone}. Single and double precision codes share the
    same parameter setting. The parallel results refer to runs with
    $\proc = 8, 64$. We avoid running the codes with 1 processor, which would have taken
    more than one day of simulation.

\begin{figure}[h!]
    \begin{minipage}{8cm}
        \begin{center}
        	\begin{tabular}{lr|rrr|}
                & \proc & CPU & $S_{\proc}$ & its (inner)\\
                single & 8 & 12856 & -- & 14 (56) \\
                single & 64 &  2698& 4.77$^{*}$ & 14 (56) \\\hline
                double & 8 & 18265  & -- & 13 (52) \\
                double & 64 & 3762  & 4.86$^{*}$ & 13 (52) \\
                \hline
                \multicolumn{5}{l}{$^{*}$ Ratio between $T_{8}$ and $T_{64}$}
            \end{tabular}
        \end{center}	
    \end{minipage}
    \hspace{-2mm}
    \begin{minipage}{6cm}
    \vspace{-2.5cm}
\includegraphics[width=7.2cm]{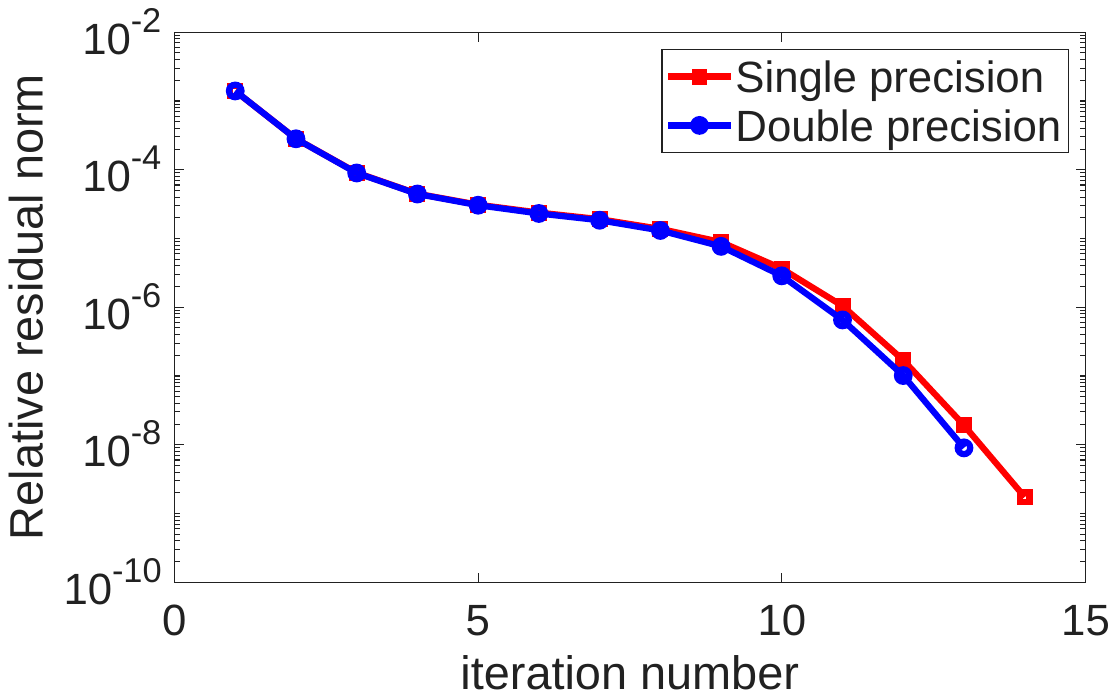}
    \vspace{-3cm}
    \end{minipage}
    \caption{Test case \#3. Simulation results (left) and convergence profiles of
    FGMRES (right).}
    \label{figconv}
\end{figure}

Employing the same parameters for both settings, the convergence profiles do not
    significantly change from double to single precision, as reported in Figure
    \ref{figconv} (right); the CPU times are 30\% lower using single precision.

Finally, in Table \ref{Single_vs_Double} we report the ratios between the CPU times required
    for the solution of the system when employing the preconditioner in double precision
    over the CPU times when employing single precision. We can observe that the mixed
    precision approach is always faster than the full double one, with peaks of about 50\%
    speed-up achieved, see, for example, test case \#2 with \proc=1. On average, we observe
    that the full double approach is 1.5 times slower than the mixed precision one.

\begin{table}[h!]
    \caption{Ratios between double and single precision CPU times, with varying number
    of processors.}
    \label{Single_vs_Double}
    \begin{center}
        \begin{tabular}{r|ccc}
            & \multicolumn{3}{c}{\proc} \\
            test case & 1 & 8 & 64 \\
            \hline
            \#1  & 1.3428  &  1.2590 &   1.1104 \\
            \#2 & 2.1484   & 1.7574  &  1.5450\\
            \#3 & -- &   1.4207 & 1.3944   \\    
        \end{tabular}
    \end{center}
\end{table}

\section{Conclusion}\label{sec_6}
In this work, we considered the numerical integration of the heat equation when
    employing a Runge--Kutta method in time. We considered a full space-time
    discretization, coupling the numerical solution and the stages for all time steps.
    The resulting system has been solved with a preconditioned iterative method, with
    the preconditioner applied in a mixed precision framework. Sequential and parallel
    numerical tests showed the robustness of the proposed mixed precision strategy, which
    is comparable to the full double precision solver in terms of number of iterations.
    Further, parallel results showed the speed-up obtained when running the iterative
    solver within a mixed precision framework.

Future work will consider the application of the mixed precision preconditioner
    to more complicated DAEs, such as differential equations describing incompressible
    fluid flow problems.

\bibliographystyle{siam}
\bibliography{references}

\end{document}